\documentclass[12pt,twoside,a4paper]{amsart}
\usepackage{amsmath,amssymb,amscd,amsfonts,latexsym,a4wide}

\newcommand{\C}{{\mathbb C}}

\newcommand{\Z}{{\mathbb Z}}
\newcommand{\Q}{{\mathbb Q}}

\renewcommand{\P}{{\mathbb P}}

\newcommand{\bx}{{\bf x}}
\newcommand{\by}{{\bf y}}

\newcommand{\bu}{{\bf u}}
\newcommand{\ba}{{\bf a}}
\newcommand{\bb}{{\bf b}}

\newcommand{\bv}{{\bf v}}
\newcommand{\bw}{{\bf w}}

\newcommand{\bA}{{\bf A}}

\newcommand{\bB}{{\bf B}}

\renewcommand{\subjclass}[1]
{\thanks{\emph{2000 Mathematics Subject Classification:}~#1}}

\renewcommand{\keywords}[1]
{\thanks{\emph{Keywords and Phrases:}~#1}}

\renewcommand{\a}{\alpha}
\newcommand{\bxi}{\xi}

\numberwithin{equation}{section}
\newcommand{\kdots}{,\ldots,}
\newcommand{\ddz}{\frac{{\rm d}}{{\rm d}z}}

\newcommand{\rank}{\mbox{$\mathrm{rank}\,$}}
\renewcommand{\det}{\mbox{$\mathrm{det}\,$}}

\newcommand{\formref}[1]{(\ref{#1})}

\newcommand{\kn}{(k^*)^n}
\newcommand{\Kn}{(K^*)^n}

\newcommand{\Ga}{\Gamma}

\renewcommand{\leq}{\leqslant}
\renewcommand{\geq}{\geqslant}

\title[Linear equations with unknowns from a multiplicative group]{Linear equations
with unknowns from a multiplicative group in a function field}


\author[J.-H.~EVERTSE]{Jan-Hendrik~EVERTSE}
\author[U.~ZANNIER]{Umberto~ZANNIER}
\address{Jan-Hendrik Evertse,\newline
\indent Universiteit Leiden, Mathematisch Instituut,\newline
\indent Postbus 9512, 2300 RA Leiden, The Netherlands}
\email{evertse@math.leidenuniv.nl}
\address{Umberto Zannier,\newline
\indent Universit\`{a} di Udine, 
Dipartimento di Matematica e Informatica,
\newline
\indent Via delle Scienze, 206, 33100 Udine, Italy}
\email{zannier@dimi.uniud.it} 

\begin{document}

\maketitle

\begin{center}
\emph{To Professor Wolfgang Schmidt on his 70th birthday}
\\[0.5cm]
\end{center}
\section{Introduction}

Let $K$ be a field of characteristic 0, and $n$ an integer $\geq 2$.
Denote by $\Kn$ the $n$-fold direct product
of the multiplicative group
$K^*$. Thus, the group operation of $\Kn$ is coordinatewise
multiplication $(x_1\kdots x_n)\cdot (y_1\kdots y_n)=
(x_1y_1\kdots x_ny_n)$. We write $(x_1\kdots x_n)^m:= (x_1^m\kdots x_n^m)$
for $m\in\Z$.
We will often denote elements of $\Kn$ by bold face characters
$\bx$, $\by$, etc.

Evertse, Schlickewei and Schmidt
\cite{EvertseSchlickeweiSchmidt02}
proved that if
$\Ga$ is a subgroup of $\Kn$ of finite rank $r$ and
$a_1\kdots a_n$ are non-zero
elements of $K$, then
the equation
\begin{equation}
\label{1.3}
a_1x_1+\cdots +a_nx_n=1\quad\mbox{in $\bx =(x_1\kdots x_n)\in\Ga$}
\end{equation}
has at most $e^{(6n)^{3n}(r+1)}$ non-degenerate solutions, i.e.,
solutions with
\begin{equation}\label{1.4}
\sum_{i\in I} a_ix_i\not=0\quad\mbox{for each proper, non-empty subset $I$
of $\{ 1\kdots n\}$.}
\end{equation}

In the present paper, we prove a function field analogue of this result.
Thus, let $k$ be an algebraically closed field of characteristic $0$
and let $K$ be a trans\-cendental field extension of $k$, where we allow
the transcendence degree to be arbitrarily large.
Let $\Ga$ be a subgroup of $\Kn$ such that $\kn\subset\Ga$ and
such that $\Ga /\kn$ has finite rank.
This means that
there are $\ba_1\kdots \ba_r\in\Ga$
such that for every $\bx\in\Ga$ there are integers $m$,
$z_1\kdots z_r$ with $m>0$ and $\bxi\in \kn$ such that
$\bx^m=\bxi\cdot \ba_1^{z_1}\cdots \ba_r^{z_r}$. If $\Ga =\kn$
then $\Ga /\kn$ has rank $0$; otherwise,
$\rank (\Ga /\kn )$ is the smallest $r$ for which $\ba_1\kdots \ba_r$
as above exist.

We deal again with equation \formref{1.3} in solutions
$(x_1\kdots x_n)\in\Ga$
with coefficients $a_1\kdots a_n\in K^*$.
We mention that in the situation we are considering now,
\formref{1.3} might have infinitely many non-degenerate solutions.
But one can show that the set of non-degenerate solutions of
\formref{1.3} is contained in finitely many $\kn$-cosets, i.e., in finitely
many sets of the shape $\bb\cdot \kn =\{ \bb\cdot \bxi :\, \bxi\in\kn\}$
with $\bb\in\Ga$. More precisely, we prove the following:
\\[0.3cm]
{\bf Theorem.}
\emph{Let $k$ be an algebraically closed field of characteristic $0$,
let $K$ be a transcendental extension of $k$, let $n\geq 2$,
let $a_1\kdots a_n\in K^*$, and let
$\Ga$ be a subgroup of $\Kn$
satisfying
\begin{equation}\label{1.1}
\kn\subset\Ga\, ,\quad \rank (\Ga /\kn )=r>0\, .
\end{equation}
Then the set of non-degenerate solutions of equation \formref{1.3}
is contained in the union of not more than}
\begin{equation}\label{1.5}
\sum_{i=2}^{n+1} {i\choose 2}^r -n+1
\end{equation}
\emph{$\kn$-cosets.}
\\[0.3cm]

We mention that Bombieri, Mueller and Zannier \cite{BombieriMuellerZannier01}
by means of a new approach gave a rather sharp upper bound for the number
of solutions of polynomial-exponential equations in
one variable over function fields.
Their approach and result were extended by
Zannier \cite{Zannier03}
to polynomial-exponential equations over function fields in several variables.
Our proof heavily uses the arguments from this last paper.

Let us consider the case $n=2$, that is, consider the equation
\begin{equation}\label{1.7}
a_1x_1+a_2x_2=1\quad\mbox{in $(x_1,x_2)\in\Ga$,}
\end{equation}
where $\Ga$, $a_1,a_2$ satisfy the hypotheses of the Theorem with $n=2$.
It is easy to check that all solutions $(x_1,x_2)$ of \formref{1.7}
with $a_1x_1/a_2x_2\in k^*$ (if any such exist)
lie in the same $(k^*)^2$-coset,
while any two different solutions $(x_1,x_2)$
with $a_1x_1/a_2x_2\not\in k^*$ lie in different $(k^*)^2$-cosets.
So our Theorem
implies that \formref{1.7} has at most $3^r$ solutions $(x_1,x_2)$
with $a_1x_1/a_2x_2\not\in k^*$. This is a slight extension
of a result by Zannier \cite{Zannier03} who obtained the same
upper bound, but for groups $\Ga =\Ga_1\times\Ga_1$ where $\Ga_1$
is a subgroup of $K^*$.

The formulation of our Theorem was inspired by Mueller
\cite{Mueller00}. She proved that
if $S$ is a finite set of places of the
rational function field $k(z)$, if $\Ga = U_S^n$ is the $n$-fold direct
product of the group of $S$-units in $k(z)^*$, and if $a_1\kdots a_n\in k(z)^*$,
then the set of non-degenerate solutions of \formref{1.1} is contained
in the union of not more than $\big( e(n+1)!/2\big)^{n(2|S|+1)}$
$\kn$-cosets.

Evertse and Gy\H{o}ry \cite{EvertseGyory88} also considered equation
\formref{1.3} with $\Ga = U_S^n$, but in the more general situation
that $S$ is a finite set of places in any finite extension $K$ of $k(z)$.
They showed that if $K$ has genus $g$
and if $a_1\kdots a_n\in K^*$
then the set of solutions $\bx\in U_S^n$
of \formref{1.3} with $(a_1x_1\kdots a_nx_n)\not\in (k^*)^n$ is
contained in the union of not more than
\[
\log (g+2)\cdot \big( e(n+1)\big)^{(n+1)|S|+2}
\]
proper linear subspaces of $K^n$.

We mention that in general $\rank U_S^n\leq n(|S|-1)$ but that in
contrast to number fields, equality need not hold.
From our Theorem we can deduce the following
result, which removes the dependence on the genus $g$,
and replaces the dependence on
$|S|$ by one on the rank.
\\[0.3cm]
{\bf Corollary.} \emph{Let $k$, $K$, $n$, $a_1\kdots a_n$, $\Ga$, $r$
be as in the Theorem.
Then the set of solutions of \formref{1.3} with
$(a_1x_1\kdots a_nx_n)\not\in (k^*)^n$ is contained in the union of
not more than
\begin{equation}\label{1.6}
\sum_{i=2}^{n+1} {i\choose 2}^r + 2^n-2n-1
\end{equation}
proper linear subspaces of $K^n$.}
\\[0.3cm]
In Section 2 we prove some auxiliary results for formal power series,
in Section 3 we prove our Theorem in the case that $K$ has transcendence
degree $1$ over $k$,
in Section 4 we extend this to the general case
that $K$ is an arbitrary transcendental extension of $k$,
and in Section 5 we deduce the Corollary.
\\[0.5cm]
\section{Results for formal power series}

Let $k$ be an algebraically closed field of characteristic $0$.
Let $z$ be an indeterminate. Denote as usual by $k[[z]]$
the ring of formal power series over $k$ and by
$k((z))$ its quotient field.
Thus, $k((z))$ consists of series $\sum_{i\geq i_0} c_iz^i$
with $i_0\in\Z$ and $c_i\in k$ for $i\geq i_0$.
We endow $k((z))$ with a derivation
$\ddz :\sum_{i\geq i_0} c_iz^i\mapsto
\sum_{i\geq i_0} ic_iz^{i-1}$.
Let $1+zk[[z]]$ denote the set of all formal power series
of the shape
$1+c_1z+c_2z^2+\cdots$ with $c_1,c_2,\ldots \in k$.
Clearly, $1+zk[[z]]$ a multiplicative group.
For $f\in 1+zk[[z]]$, $u\in k$ we define
\begin{equation}\label{2.000}
f^u := \sum_{i=0}^{\infty} {u\choose i}(f-1)^i\, ,
\end{equation}
where ${u\choose 0}=1$ and
${u\choose i}=u(u-1)\cdots (u-i+1)/i!$ for $i>0$.
Thus, $f^u$ is a well-defined element of $1+zk[[z]]$.
This definition of $f^u$ coincides
with the usual one for $u=0,1,2,\ldots$.
We have $\ddz f^u =uf^{u-1}\ddz f$ and moreover,
\begin{equation}\label{2.0}
\left\{
\begin{array}{l}
\mbox{$(fg)^u=f^ug^u$ for $f,g\in 1+zk[[z]]$, $u\in k$;}\\[0.1cm]
\mbox{$f^{u+v}=f^uf^v$ and
$(f^u)^v=f^{uv}$ for $f\in 1+zk[[z]]$, $u,v\in k$.}
\end{array}\right.
\end{equation}
(One may verify  \formref{2.0} by taking logarithmic derivatives
and using that two series in $1+zk[[z]]$ are equal if and only
if their logarithmic derivatives are equal).
We endow $(1+zk[[z]])^r$ with the usual coordinatewise multiplication.
Given $\bB =(b_1\kdots b_r)\in (1+zk[[z]])^r$,
we define $\bB^u:= (b_1^u\kdots b_r^u)$ for $u\in k$ and
$\bB^{\bu}:= b_1^{u_1}\cdots b_r^{u_r}$ for $\bu =(u_1\kdots u_r)\in k^r$.
Thus, $\bB^u\in (1+zk[[z]])^r$ and $\bB^{\bu}\in 1+zk[[z]]$.

Let $h,r$ be integers with $h\geq 2$, $r\geq 1$.
Further, let $a_1\kdots a_h$ be elements of $k[[z]]$
which are algebraic over
the field of rational functions
$k(z)$ and which are not divisible by $z$,
and let $\a_{ij}$ $(i=1\kdots h,\, j=1\kdots r)$
be elements of $1+zk[[z]]$ which are algebraic over $k(z)$.
Put $\bA_i :=(\a_{i1}\kdots \a_{ir})$ ($i=1\kdots h$).
Define
\[
R:=\{ \bu\in k^r:\,\mbox{$a_1\bA_1^{\bu}\kdots
a_h\bA_h^{\bu}$ are linearly dependent over $k$.}\}
\]
By a \emph{class} we mean a set $R'\subset k^r$ with the property that
there are a subset $J$ of $\{ 1\kdots h\}$ and $\bu_0\in\Q^r$
such that for every $\bu\in R'$ the following holds:
\begin{equation}\label{2.00}
\left\{
\begin{array}{l}
\mbox{$a_i\bA_i^{\bu}$ $(i\in J)$ are linearly dependent over $k$;}\\
\mbox{$\big( \bA_i\bA_j^{-1}\big)^{\bu -\bu_0}=1$ for all $i,j\in J$.}
\end{array}\right.
\end{equation}
\\[0.15cm]
{\bf Lemma 1.}
\emph{$R$ is the union of finitely many classes.}
\\[0.3cm]
{\bf Proof.}
This is basically a special case of \cite[Lemma 1]{Zannier03}. In the proof
of that lemma, it was assumed that $k=\C$, and that the $a_i$ and
$\a_{ij}$ are holomorphic functions in the variable $z$
which are algebraic
over $\C (z)$
and which are
defined and have no zeros
on a simply connected open subset $\Omega$ of $\C$.
It was shown that
provided $k=\C$, this was no loss of generality. The argument remains
precisely the same if one allows $k$ to be an arbitrary algebraically closed
field of characteristic $0$ and if one takes for the $a_i$
power series from $k[[z]]$ which are algebraic over $k(z)$ and which are
not divisible by $z$,
and for the $\a_{ij}$
power series from $1+zk[[z]]$ which are algebraic over $k(z)$.

We mention that in \cite{Zannier03} the definition of a class
is slightly different from \formref{2.00}, allowing
$\big( \bA_i\bA_j^{-1}\big)^{\bu -\bu_0}\in k^*$ for all $i,j\in J$.
But in our situation this implies automatically that
$\big( \bA_i\bA_j^{-1}\big)^{\bu -\bu_0}=1$
since $\big( \bA_i\bA_j^{-1}\big)^{\bu -\bu_0}\in 1+zk[[z]]$.
\qed
\\[0.3cm]

We now impose some further restriction on the $\a_{ij}$ and prove
a more precise result.
Namely, we assume that
\begin{equation}\label{2.1}
\{ \bu\in k^r:\, (\bA_i\cdot\bA_h^{-1})^{\bu}=1\,\,\,
\mbox{for $i=1\kdots h$}\}=\{{\bf 0}\}\, .
\end{equation}
Let $S$ be the set of $\bu\in k^r$ such that there are
$\xi_1\kdots \xi_h\in k$
with
\begin{eqnarray}\label{2.2}
&&\sum_{i=1}^h \xi_ia_i\bA_i^{\bf u}=0,\,
\\
\label{2.3}
&&\sum_{i\in I} \xi_ia_i\bA_i^{\bf u}\not= 0
\quad\mbox{for each proper, non-empty subset $I$ of $\{ 1\kdots h\}$.}
\end{eqnarray}
\\[0.15cm]
{\bf Lemma 2.} \emph{Assume \formref{2.1}. Then $S$ is finite.}
\\[0.3cm]
{\bf Proof.} We prove a slightly stronger statement. We partition
$\{ 1\kdots h\}$ into subsets $I_1\kdots I_s$ such that $\bA_i=\bA_j$
if and only if $i,j$ belong to the same set $I_l$ for some
$l\in\{1\kdots s\}$.
Let $\tilde{S}$ be the set of $\bu\in k^r$ satisfying \formref{2.2}
and, instead of \formref{2.3},
\begin{equation}\label{2.4}
\sum_{i\in I} \xi_ia_i\bA_i^{\bf u}\not= 0
\end{equation}
for each proper, non-empty subset $I$ of $\{ 1\kdots h\}$ which is a union
of some of the sets $I_1\kdots I_s$. We prove that $\tilde{S}$ is finite.
This clearly suffices.

We proceed by induction on $p:=h+s$. Notice that from assumption \formref{2.1}
it follows that $h\geq 2$ and $s\geq 2$. First let $h=2$, $s=2$, i.e.,
$p=4$. Thus, $\tilde{S}$ is the set of $\bu\in k^r$ for which there are
non-zero $\xi_1,\xi_2\in k$ with
$\xi_1a_1\bA_1^{\bf u}+\xi_2a_2\bA_2^{\bf u}=0$.
Then for $\bu\in \tilde{S}$ we have
\[
\big(\bA_1\cdot \bA_2^{-1}\big)^{\bf u}=\xi (a_2a_1^{-1})
\]
with $\xi \in k^*$. Consequently,
$\big(\bA_1\cdot \bA_2^{-1}\big)^{\bu_2-\bu_1}\in k^*$ for any
$\bu_1,\bu_2\in\tilde{S}$. But then for $\bu_1,\bu_2\in\tilde{S}$ 
we must have 
$\big(\bA_1\cdot \bA_2^{-1}\big)^{\bu_2-\bu_1}=1$
since 
$\big(\bA_1\cdot \bA_2^{-1}\big)^{\bu_2-\bu_1}\in 1+zk[[z]]$.
In view of assumption \formref{2.1} this implies that
$\tilde{S}$ consists of at most one element.

Now let $p>4$ and assume Lemma 1 is true for all pairs $(h,s)$ with
$h\geq 2$, $s\geq 2$ and $h+s<p$. We apply Lemma 1 above.
Clearly, $S$ is contained in the set $R$ defined above, and therefore,
$\tilde{S}$ is the union of finitely many sets $\tilde{S}\cap R'$
where $R'$ is a class as defined above.
So we have to show that each such set $\tilde{S}\cap R'$ is finite.

Thus let $S':= \tilde{S}\cap R'$, where $R'$ is a class as above.
Let $J$ be the corresponding subset of $\{ 1\kdots h\}$,
and $\bu_0\in\Q^r$ the corresponding vector, such that \formref{2.00} holds.
We distinguish two cases. First suppose that $J$ is contained in some
set $I_l$. Then the elements $a_j$ $(j\in J)$ are linearly
dependent over $k$. There is a proper subset $J'$ of $J$ such that
$a_j$ $(j\in J')$ are linearly independent over $k$ and such that
each $a_j$ with $j\in J\backslash J'$ can be expressed as a linear combination
over $k$ of the $a_j$ with $j\in J'$. By substituting these linear
combinations into \formref{2.2}, \formref{2.4},
we obtain similar conditions, but with $I_l$ replaced by the smaller
set obtained by removing from $I_l$ the elements from $J\backslash J'$.
This reduces the value of the number of terms $h$.
Further, condition \formref{2.1} remains valid.
Thus we may apply the induction hypothesis, and conclude that $S'$ is finite.

Now assume that $J$ is not contained in one of the sets $I_l$.
We transform our present situation into a new one with 
instead of $I_1\kdots I_s$ a partition of $\{ 1\kdots h\}$
into fewer than $s$ sets.
Then again, the induction hypothesis is applicable.

There are $i,j\in J$ with $\bA_i\not=\bA_j$,
say $i\in I_{l_1}$ and $j\in I_{l_2}$.
Further, there is $\bu_0\in \Q^r$
such that $\big( \bA_i\bA_j^{-1}\big)^{\bu-\bu_0}=1$ for $\bu\in S'$.
According to an argument in the proof of Lemma 1 of \cite{Zannier03},
the set of $\bu\in k^r$ with $\big( \bA_i\bA_j^{-1}\big)^{\bu}=1$
is a linear subspace $V$ of $ k^r$ which is defined over $\Q$.
Let $\bv_1\kdots \bv_{r'}$ be a basis of $V$ contained in $\Z^r$.
Thus, each $\bu\in S'$ can be expressed uniquely as
\begin{equation}\label{2.10}
\bu_0 +w_1\bv_1+\cdots +w_{r'}\bv_{r'}\quad
\mbox{with $\bw=(w_1\kdots w_{r'})\in k^{r'}$.}
\end{equation}
Now define
\[
b_q := a_q\bA_q^{\bu_0}\, ,\quad
\bB_q :=(\bA_q^{\bv_1}\kdots \bA_q^{\bv_{r'}})\quad
(q=1\kdots h)\, .
\]
Thus, for $\bu\in S'$ we have
\begin{equation}\label{2.11}
a_q\bA_q^{\bu}= b_q\bB_q^{\bw}\quad\mbox{for $q=1\kdots h$.}
\end{equation}
Clearly, $b_q\in k[[z]]$ and the coordinates of $\bB_q$
belong to $1+zk[[z]]$, for $q=1\kdots h$.
Further, $b_q$,
and the coordinates of $\bB_q$ $(q=1\kdots h)$ are algebraic over $k(z)$
since $\bu_0\in\Q^r$ and since $\bv_1\kdots\bv_{r'}\in\Z^r$.

From
the definition of $\bB_q$ $(q=1\kdots h)$ it follows that
if $(\bB_q\bB_h^{-1})^{\bw}=1$ for $q=1\kdots h$, then
$(\bA_q\bA_h^{-1})^{\sum_j w_j\bv_j}=1$ for $q=1\kdots h$, which by
\formref{2.1} implies $\sum_j w_j\bv_j={\bf 0}$ and so $\bw ={\bf 0}$.
Therefore, condition \formref{2.1} remains valid if we replace
$\bA_q$ by $\bB_q$ for $q=1\kdots h$.

It is important to notice that $\bB_{q_1}=\bB_{q_2}$
for any $q_1,q_2\in I_{l_1}\cup I_{l_2}$.
Further, for each $l\not= l_1,l_2$, we have that
$\bB_{q_1}=\bB_{q_2}$
for any $q_1,q_2\in I_l$.

Lastly, if $\bu\in S'$ then by substituting \formref{2.11}
into \formref{2.2}, \formref{2.4}, we obtain
that there are $\xi_1\kdots\xi_h\in k^*$
such that
$\sum_{q=1}^h \xi_q b_q\bB_q^{\bw}=0$ and
$\sum_{q\in I} \xi_q b_q\bB_q^{\bw}\not=0$
for each proper subset $I$ of $\{ 1\kdots h\}$ which is a union of some of
the sets from $I_{l_1}\cup I_{l_2}$,
$I_l$ $(l=1\kdots s,\, l\not= l_1,l_2)$.
Thus, each $\bu\in S'$ corresponds by means of \formref{2.10}
to $\bw\in k^{r'}$ which satisfies similar conditions as $\bu$,
but with instead of $I_1\kdots I_s$ a partition of $\{ 1\kdots h\}$
into $s-1$ sets.
Now by the induction
hypothesis, the set of $\bw$ is finite, and therefore, $S'$ is finite.
This proves Lemma 2. \qed
\\[0.4cm]

We now proceed to estimate the cardinality of $S$. We need a few
auxiliary results. For any subset $A$ of $k[[z]]$,
we denote by $\rank_k A$ the cardinality of
a maximal  $k$-linearly independent subset of $A$.
For each subset $I$ of $\{ 1\kdots h\}$ and each integer $t$ with
$1\leq t\leq h-1$, we define the set
\begin{equation}\label{2.5}
V(I,t)=\{ \bu\in k^r\, : \rank_k\{ a_i\bA_i^{\bu}:\, i\in I\}\leq t\}\, .
\end{equation}
Clearly, $V(I,t)=k^r$ if $t\geq |I|$.
\\[0.3cm]
{\bf Lemma 3.} \emph{Let $I,t$ be as above and assume that $t<|I|$. Then
$V(I,t)$ is the set of common zeros in $k^r$ of a system
of polynomials in $k[X_1\kdots X_r]$, each of total degree at most
${t+1\choose 2}$.}
\\[0.3cm]
{\bf Proof.}
The vector $\bu$ belongs to $V(I,t)$ if and only if each $t+1$-tuple
among the functions $a_i\bA_i^{\bu}$ $(i\in I)$ is linearly
dependent over $k$, that is,
if and only if for each subset $J=\{ i_0\kdots i_t\}$ of $I$
of cardinality $t+1$, the Wronskian determinant
\[
\det\left( \big(\ddz\big)^i
a_{i_j}\bA_{i_j}^{\bu}\right)_{i,j=0\kdots t}
\]
is identically $0$ as a function of $z$.
By an argument completely similar to that in the proof
of Proposition 1 of \cite{Zannier03}, one shows that the latter condition
is equivalent to $\bu$ being a common zero of some finite set of
polynomials of degree $\leq {t+1\choose 2}$. This proves Lemma 3.\qed
\\[0.3cm]
{\bf Lemma 4.} \emph{
$\bu\in S$ if and only if
\begin{eqnarray}
\label{2.7}
&&\rank_k\{ a_i\bA_i^{\bu}:\, i\in I\}+
\rank_k\{ a_i\bA_i^{\bu}:\, i\not\in I\}\\
&&\qquad\qquad\qquad\qquad\qquad\qquad >
\rank_k\{ a_i\bA_i^{\bu}:\, i=1\kdots h\}
\nonumber
\end{eqnarray}
for each proper, non-empty subset $I$ of $\{ 1\kdots h\}$.}
\\[0.3cm]
{\bf Proof.} First let $\bu\in S$.
Take a proper, non-empty subset $I$ of
$\{ 1\kdots h\}$. From \formref{2.2}, \formref{2.3} it follows that
there are $\xi_1\kdots \xi_h\in k$ such that
\[
\sum_{i\in I} \xi_i a_i\bA_i^{\bu} = -\sum_{i\not\in I} \xi_ia_i\bA_i^{\bu}
\not= 0
\]
and therefore the $k$-vector spaces spanned by $\{ a_i\bA_i^{\bu}:\, i\in I\}$,
$\{ a_i\bA_i^{\bu}:\, i\not\in I\}$, respectively,
have non-trivial intersection. This implies \formref{2.7}.

Now let $\bu\in k^r$ be such that \formref{2.7} holds
for every proper, non-empty subset $I$ of $\{ 1\kdots h\}$.
Let $W$ be the
vector space of $\bxi =(\xi_1\kdots \xi_h)\in k^h$ with
$\sum_{i=1}^h \xi_i a_i\bA_i^{\bu}=0$. Further, for a proper, non-empty subset
$I$ of $\{ 1\kdots h\}$, let $W(I)$ be the vector space of
$\bxi =(\xi_1\kdots \xi_h)\in k^h$ with $\sum_{i\in I} \xi_i a_i\bA_i^{\bu}=0$
and $\sum_{i\not\in I} \xi_i a_i\bA_i^{\bu}=0$.
Given a proper, non-empty subset $I$ of $\{ 1\kdots h\}$,
it follows from \formref{2.7} that there are $\xi_1\kdots\xi_h\in k$
with $\sum_{i\in I} \xi_i a_i\bA_i^{\bu} =
-\sum_{i\not\in I} \xi_i a_i\bA_i^{\bu}\not= 0$; hence $W(I)$ is a proper
linear subspace of $W$. It follows that there is $\bxi\in W$ with
$\bxi\not\in W(I)$ for each proper, non-empty subset $I$ of
$\{ 1\kdots h\}$. This means precisely that $\bu\in S$.\qed
\\[0.3cm]
{\bf Proposition.}
\emph{Assume \formref{2.1}. Then $|S|\leq \sum_{p=2}^h {p\choose 2}^r -h+2$.}
\\[0.3cm]
{\bf Proof.} For $t=1\kdots h-1$, let $T_t=V(\{ 1\kdots h\}, t)$
(that is the set of $\bu\in k^r$ 
with $\rank_k \{ a_i\bA_i^{\bu}:\, i=1\kdots h\}\leq t$)
and let $S_t$ be the set of $\bu\in S$ such that
$\rank_k\{ a_i\bA_i^{\bu}:\, i=1\kdots h\}=t$.
By \formref{2.7},
$\rank_k\{ a_i\bA_i^{\bu}:\, i=1\kdots h\}<h$, so
$S=S_1\cup\cdots\cup S_{h-1}$.
We show by induction on $t=1\kdots h-1$ that
\begin{equation}
\label{2.8}
|S_1\cup\cdots\cup S_t|\leq \sum_{p=1}^t {p+1\choose 2}^r -t+1\, .
\end{equation}
Taking $t=h-1$, our Proposition follows.

First let $t=1$. Let $\bu_1,\bu_2\in S_1$.
Then $(a_ia_h^{-1})(\bA_i\bA_h^{-1})^{\bu_j}\in k^*$ for $i=1\kdots h$, $j=1,2$,
which implies
$(\bA_i\bA_h^{-1})^{\bu_1-\bu_2}\in k^*$ for $i=1\kdots h$.
But then $(\bA_i\bA_h^{-1})^{\bu_1-\bu_2}=1$
since $(\bA_i\bA_h^{-1})^{\bu_1-\bu_2}\in 1+zk[[z]]$
for $i=1\kdots h$.
Now assumption \formref{2.1} gives
$\bu_1=\bu_2$.
So $|S_1|=1$ which implies \formref{2.8} for $t=1$.

Now assume that $2\leq t\leq h-1$ and that \formref{2.8} is true
with $t$ replaced by any number $t'$ with $1\leq t'<t$.
By Lemma 3, $T_t$ is an algebraic subvariety of $k^r$,
being the set of common zeros of a system of polynomials
of degree not exceeding ${t+1\choose 2}$. By the last part of the
proof of Proposition 1 of \cite{Zannier03}, $T_t$ has at most
${t+1\choose 2}^r$ irreducible components.

We first show that $T_t\backslash S_t$ is a finite union of proper algebraic
subvarieties of $T_t$.
Notice that $\bu\in T_t\backslash S_t$ if and only if either
$\rank_k\{ a_i\bA_i^{\bu}:\, i=1\kdots h\}\leq t-1$
or (by Lemma 4) there are a proper, non-empty subset $I$ of $\{ 1\kdots h\}$
and an integer $q$ with $1\leq q\leq t-1$ such that
$\rank_k\{ a_i\bA_i^{\bu}:\, i\in I\}\leq q$ and
$\rank_k\{ a_i\bA_i^{\bu}:\, i\not\in I\}\leq t-q$.
This means that $T_t\backslash S_t$ is equal to the union
of $T_{t-1}$ and of all sets $V(I,q)\cap V(\{ 1\kdots h\}\backslash I,t-q)$
with $I$ running through the proper, non-empty subsets
of $\{ 1\kdots h\}$ and $q$ running through the integers
with $1\leq q\leq t-1$. By Lemma 3
these sets are all subvarieties of $T_t$.

Now by Lemma 2 $S_t$ is finite,
hence each element of $S_t$ is an irreducible component (in fact
an isolated point) of $T_t$.
So $|S_t|\leq {t+1\choose 2}^r$.
Now two cases may occur.

If $T_t=S_t$ then $S_{t'}=\emptyset$ for $t'=1\kdots t-1$ and so
$|S_1\cup\cdots\cup S_t|=|S_t|\leq {t+1\choose 2}^r$. This certainly
implies \formref{2.8}.

If $S_t$ is strictly smaller than $T_t$ then $T_t\backslash S_t$ has at least
one irreducible component. But then $|S_t|\leq {t+1\choose 2}^r-1$.
In conjunction with the induction hypothesis this gives
\begin{eqnarray*}
|S_1\cup\cdots\cup S_t| &=& |S_1\cup\cdots\cup S_{t-1}|+|S_t|\\
&\leq&
\sum_{p=1}^{t-1}{p+1\choose 2}^r-t+2 +{t+1\choose 2}^r-1
\end{eqnarray*}
which implies again \formref{2.8}.

This completes the proof of our induction step, hence of our Proposition.\qed
\\[0.5cm]
\section{Proof of the Theorem for transcendence degree $1$}

We prove the Theorem in the special case that $K$ has transcendence degree
$1$ over $k$.
For convenience we put $N:=\sum_{i=2}^{n+1} {i\choose 2}^r-n+1$.

We start with some reductions.
There are $\ba_j=(\a_{1j}\kdots\a_{nj})\in \Ga$
$(j=1\kdots r)$ such that for each $\bx\in\Ga$ there are integers
$m,\, w_1\kdots w_r$ with $m>0$, and $\xi =(\xi_1\kdots \xi_n)\in\kn$
such that $\bx^m = \xi\cdot\ba_1^{w_1}\cdots\ba_r^{w_r}$.
Let $L$ be the extension of $k$ generated by $a_1\kdots a_n$
and the $\a_{ij}$ ($i=1\kdots n$, $j=1\kdots r$).
Then $L$ is the function field of
a smooth projective algebraic curve $C$ defined over $k$.
Choose $z\in L$, $z\not\in k$, such that the map
$z: C\to \P_1(k)=k\cup\{\infty\}$ is unramified at $0$
and such that none of the functions $a_i$, $\a_{ij}$ has a zero or pole
in any of the points from $z^{-1}(0)$.
Thus, $L$ can be embedded into $k((z))$, and the $a_i$ and $\a_{ij}$
may be viewed as elements of $k[[z]]$ not divisible by $z$.
By multiplying the $\a_{ij}$ with appropriate constants from $k^*$,
which we are free to do, we may assume without loss of generality
that the $\a_{ij}$ belong to $1+zk[[z]]$.

Making the asumptions for the $a_i$ and $\a_{ij}$ just mentioned,
we can apply our Proposition.
The functions $\a_{ij}^u$
$(u\in k)$ are defined uniquely by means of \formref{2.000}.
Therefore, we can express each $\bx\in\Ga$ as
$\xi\cdot \ba_1^{u_1}\cdots\ba_r^{u_r}$ with $u_1\kdots u_r\in\Q$
and with $\xi =(\xi_1\kdots\xi_n)\in \kn$.
Putting
$\bA_i:= (\a_{i1}\kdots \a_{ir})$ $(i=1\kdots n)$,
we can rewrite this as
\begin{equation}\label{3.0}
\bx =(\xi_1\bA_1^{\bu}\kdots \xi_n\bA_n^{\bu})
\end{equation}
with $\xi_1\kdots\xi_n\in k^*$, $\bu =(u_1\kdots u_r)\in\Q^r$.
Putting in addition $h:=n+1$,  $\bA_h:= (1\kdots 1)$ ($r$ times $1$),
$a_h:=-1$, $\xi_h:= 1$ we obtain that
if $\bx\in\Ga$ is a
non-degenerate solution of \formref{1.3} then
\begin{eqnarray}
\label{3.1}
&&\sum_{i=1}^h \xi_i a_i\bA_i^{\bu}=0\, ,
\\[0.1cm]
\label{3.2}
&&\sum_{i\in I} \xi_i a_i\bA_i^{\bu}\not=0
\quad\mbox{for each proper, non-empty subset $I$ of $\{ 1\kdots h\}$.}
\end{eqnarray}
It remains to verify condition \formref{2.1}.
According to an argument in the proof of
Lemma 1 of \cite{Zannier03}, the set of $\bu\in k^r$ such that
$(\bA_i\bA_h^{-1})^{\bu}=1$ for $i=1\kdots h$ is a linear subspace
of $k^r$, say $V$, which is defined over $\Q$.
Now if $\bu =(u_1\kdots u_r)\in V\cap\Q^r$,
then $\bA_i^{\bu}=1$ for $i=1\kdots n$ since $\bA_h =(1\kdots 1)$,
and therefore
$\ba_1^{u_1}\cdots \ba_r^{u_r}=(1\kdots 1)$.
This implies $\bu ={\bf 0}$, since otherwise $\rank (\Ga /\kn )$
would have been smaller than $r$.
Hence $V\cap\Q^r =\{ {\bf 0}\}$ and therefore, $V=\{ {\bf 0}\}$ since $V$
is defined over $\Q$. This implies \formref{2.1}.

As observed above, if $\bx\in\Ga$ is a non-degenerate solution of \formref{1.3},
then $\bu$ satisfies \formref{3.1},\formref{3.2}, which means that $\bu$
belongs to the set $S$ given by \formref{2.2}, \formref{2.3}.
So by the Proposition, we have at most
$N$ possibilities for $\bu$. Then according to \formref{3.0},
the non-degenerate solutions $\bx$ of \formref{1.3}
lie in at most $N$ $\kn$-cosets.
This completes the proof of our Theorem
in the special case that $K$ has transcendence degree $1$ over $k$.\qed
\\[0.5cm]
\section{Proof of the Theorem in the general case}

We prove our Theorem in the general case, i.e.,
that the field $K$ is an arbitrary
transcendental extension of $k$.
As before, we denote $N:=\sum_{i=2}^{n+1}{i\choose 2}^r-n+1$.

There is of course no loss of generality to assume that
$K$ is generated by the coefficients $a_1\kdots a_n$
and the coordinates of all elements of $\Ga$.
Since $\Ga$ is assumed to have rank $r$,
there are $\ba_1\kdots \ba_r\in\Ga$
such that for every $\bx\in\Ga$ there are integers $m$,
$z_1\kdots z_r$ with $m>0$ and $\bxi\in \kn$ such that
$\bx^m=\bxi\cdot \ba_1^{z_1}\cdots \ba_r^{z_r}$.
Hence $K$ is algebraic over the extension of $k$ generated by $a_1\kdots a_n$
and the coordinates of $\ba_1\kdots\ba_r$. Therefore, $K$ has finite
transcendence degree over $k$.
We will prove by induction on $d:={\rm trdeg}(K/k)$ that
for any group $\Ga$ with $\rank (\Ga /\kn )\leq r$,
the non-degenerate solutions $\bx\in\Ga$ of \formref{1.1}
lie in not more than $N$ $\kn$-cosets.
The case $d=0$ is trivial since in that case $\Ga =\kn$ and all solutions
lie in a single $\kn$-coset.
Further, the case $d=1$ has been taken care of in the previous section.
So we assume $d>1$ and that the above assertion is true up to $d-1$.

We assume by contradiction that \formref{1.3} has at least $N+1$
non-degenerate solutions, denoted $\bx_1\kdots\bx_{N+1}\in \Gamma$,
falling into pairwise distinct $\kn$-cosets.
For each such solution
$\bx_j =:(x_{1j},\ldots ,x_{nj})$ and for each  nonempty subset $I$ of
$\{1\kdots n\}$
let us consider the corresponding subsum $\sum_{i\in I}a_ix_{ij}$,
which we
denote $\sigma_{(j,I)}$. In this way we obtain finitely many elements
$\sigma_{(j,I)}\in
K$, none of which vanishes, since the solutions are non-degenerate.

Further, let $\bx_u,\bx_v$ be distinct solutions,
with $1\leq u\not= v\leq N+1$.
Since the solutions lie in distinct $(k^*)^n$-cosets, for some
$i\in\{1,\ldots ,n\}$ the ratio
$x_{iu}/x_{iv}$ does not lie in
$k$.  For each pair $(u,v)$ as
above let us pick one such index $i=i(u,v)$ and let us put
$\tau_{(u,v)} := x_{iu}/x_{iv}\in K^*\setminus k^*$.

We are going to ``specialize" such elements of $K$, getting corresponding
elements of a field with smaller transcendence degree and obtaining
eventually a
contradiction. We shall formulate the specialization argument in geometric
terms.

Let $\tilde{K}$ be the extension of $k$ generated by $a_1\kdots a_n$
and by the coordinates of $\bx_1\kdots\bx_{N+1}$.
Thus $\tilde{K}$ is finitely generated over $k$.
Further, let $\tilde{\Ga}$ be the group containing $\kn$
and generated over it by $\bx_1\kdots\bx_{N+1}$.
Then $\tilde{\Ga}$ is a subgroup of $\Ga\cap (\tilde{K}^*)^n$,
and so $\rank (\tilde{\Ga})\leq r$.
Now \formref{1.3} has at least $N+1$ non-degenerate
solutions in $\tilde{\Ga}$ lying in different $\kn$-cosets.
By the induction hypothesis
this is impossible if
${\rm trdeg} (\tilde{K}/k)<d$. So ${\rm trdeg} (\tilde{K}/k)=d$.

The finitely generated extension $\tilde{K}/k$ may be viewed as the function field
of an
irreducible affine algebraic variety
$V$ over
$k$, with $d=\dim V$.
Then,  each element of
$\tilde{K}$ represents a rational function on
$V$. Let us consider irreducible closed subvarieties $W$ of $V$, with
function field
denoted $L:=k(W)$, with the following properties:
\\[0.2cm]
(A) $\dim W= d-1$.
\\
(B) There exists a point $P\in W(k)$ such that each of the (finitely many)
elements
$a_i$, $x_{ij}$ and $\sigma_{(j,I)}$,
$\tau_{(u,v)}$ constructed above is defined and nonzero at $P$;  so the
elements induce by restriction nonzero rational functions $a_i'$, $x_{ij}'$,
$\sigma'_{(j,I)}$ and $\tau'_{(u,v)}$ in $L^*=k(W)^*$;
\\
(C) None of the elements $\tau'_{(u,v)}$ lies in $k^*$.
\\[0.2cm]
We shall construct $W$ as an irreducible component of a suitable hyperplane
section of $V$.

To start with,
(A) follows from the well-known fact that any irreducible
component $W$ of any hyperplane section of $V$ has dimension $d-1$.

Let us analyze (B). Each of the elements of $\tilde{K}^*$ mentioned in
(B) may be expressed
as a
ratio of nonzero polynomials in the affine coordinates  of $V$; since these
elements
are defined and nonzero  by assumption,  none of these polynomials vanishes
identically on
$V$, so each such polynomial defines in $V$ a proper (possibly reducible)
closed subvariety.
Take now  a point $P\in V(k)$ outside the union  of these
finitely many proper subvarieties. For (B) to be verified it then plainly
suffices
that $W$ contains $P$.

Finally, let us look at (C). For each $u,v\in\{1,\ldots ,n\}$, $u\not= v$,
let $Z(u,v)$ be the variety defined in $V$ by the equation
$\tau_{(u,v)}=\tau_{(u,v)}(P)$. Since
$\tau_{(u,v)}$ is not  constant on $V$, each component of $Z(u,v)$ is a
subvariety
of  $V$  of dimension
$d-1$. Choose now
$W$ as an irreducible component through $P$ of the intersection of $V$ with a
hyperplane
$\pi$ going through
$P$, such that $W$ is not contained in any of the finitely many $Z(u,v)$. It
suffices e.g. that the hyperplane $\pi$ does not contain any irreducible
component
of any
$Z(u,v)$ and there are plenty of choices for that. (E.g. for each of the
relevant
finitely many varieties, each of dimension
$d-1\geq 1$, take a point $Q\neq P$  in it and let $\pi$ be a hyperplane
through $P$
and not containing any of the $Q$'s. Note that here we use that
$d\geq 2$.)  Since
$P\in W(k)$ and
$\tau_{(u,v)}$ is not constantly equal to
$\tau_{(u,v)}(P)$ on all of $W$ by contruction, the restriction
$\tau'_{(u,v)}$ is
not constant, as required.

Consider now the elements $\bx_j':=(x_{1j}',\ldots ,x_{nj}')\in L^n$,
$j=1,\ldots ,N+1$, where the dash denotes, as before, the restriction to
$W$ (which
by (B) is well-defined for all the functions in question).
Notice that the restriction to
$W$ is a homomorphism from the local ring of $V$ at $P$ to the local ring
of $W$ at $P$ which is contained in $L$.
This homomorphism maps $\tilde{\Ga}$ to the group
$\Ga '$ containing $(k^*)^n$, generated over it by the elements
$\bx_1'\kdots \bx_{N+1}'$.
Thus, $a_1'\kdots a_n'$ and the coordinates of the elements
from $\Ga '$ lie in $L$.
Further, $\Ga '$ is a homomorphic image of $\tilde{\Ga}$
which was in turn a subgroup of $\Ga$; therefore $\rank (\Ga '/\kn )\leq r$.
Since the
$\bx_j$ are solutions of \formref{1.3} in $\tilde{\Ga}$, the elements $\bx_j'$
are solutions of
$a_1'x_1+\cdots +a_n'x_n=1$ in $\Ga '$.
Again by (B), we have that none of the
(nonempty) subsums
$\sigma'_{(j,I)}=\sum_{i\in I}a'_ix'_{ij}$  vanishes, so these solutions are
non-degenerate.
Finally, by (C), no two solutions $\bx_u',\bx_v'$,
$1\leq u\not= v\leq N+1$,
lie in a same
$(k^*)^n$-coset of
$(L^*)^n$.
Since by (A) the field $L$ has transcendence degree $d-1$ over
$k$, this
contradicts the inductive assumption, concluding the induction step
and the proof.\qed
\\[0.5cm]
\section{Proof of the Corollary}

We keep the notation and assumptions from Section 1.
We consider the non-degenerate solutions $(x_1\kdots x_n)\in\Ga$
of \formref{1.3} such that
\begin{equation}\label{4.1}
(a_1x_1\kdots a_nx_n)\not\in (k^*)^n\, .
\end{equation}
We first show that each $\kn$-coset of such solutions
is contained in a proper linear subspace of $K^n$.
Fix a non-degenerate solution $\bx =(x_1\kdots x_n)$
of \formref{1.3} with
\formref{4.1}.
Any other solution
of \formref{1.3} in the same $\kn$-coset as $\bx$
can be expressed as $\bx\cdot \bxi =(x_1\xi_1\kdots x_n\xi_n)$
with $\bxi = (\xi_1\kdots \xi_n)\in\kn$
and $a_1x_1\xi_1+\cdots +a_nx_n\xi_n=1$.
Now the points $\bxi\in k^n$ satisfying the
latter equation lie in a proper linear subspace of $k^n$,
since otherwise 
$(a_1x_1\kdots a_nx_n)$ would be the unique
solution of a system of $n$ linearly independent linear equations 
with coefficients from $k$, 
hence $a_1x_1\kdots a_nx_n\in k$, violating \formref{4.1}.
But this implies that indeed the $\kn$-coset
$\{ \bx\cdot \bxi :\, \bxi\in\kn\}$
is contained in a proper linear subspace of $K^n$.

Now our Theorem implies that the non-degenerate solutions of \formref{1.3}
with \formref{4.1} lie in at most $\sum_{i=2}^{n+1}{i\choose 2}^r-n+1$
proper linear subspaces of $K^n$. Further,  
the degenerate
solutions of \formref{1.3} lie in at most $2^n-n-2$ proper linear
subspaces of $K^n$, each given by $\sum_{i\in I} a_ix_i=0$, where
$I$ is a subset of $\{ 1\kdots n\}$ of cardinality $\not= 0,1,n$.
By adding these two bounds our Corollary follows.\qed
\\[0.5cm]

\end{document}